\PassOptionsToPackage{hyphens}{url}
\documentclass[submission,copyright,creativecommons]{eptcs}
\providecommand{\event}{GASCOM 2026}

\usepackage[hyphenbreaks]{breakurl}
\usepackage[utf8]{inputenc}
\usepackage[T1]{fontenc}
\usepackage[small]{caption}
\usepackage{booktabs}
\usepackage{graphicx}
\usepackage{longtable}
\usepackage{makecell}
\usepackage{amsmath}
\usepackage{amsfonts}
\usepackage{gensymb}
\usepackage{hyperref}
\usepackage{float}
\usepackage{subcaption}
\usepackage{ulem}
\usepackage{tikz}

\usepackage{todonotes}

\tikzstyle{arrow} = [ultra thick,->,>=stealth]
\tikzstyle{arrowthick} = [decoration={markings,mark=at position 1 with {\arrow[scale=3,>=stealth]{>}}},postaction={decorate}]

\title{Generation of Maximal Snake Polyominoes  Using a Deep Neural Network}
\author{Benjamin Gauthier \hfill Alain Goupil \hfill Fadel Toure
\institute{
Université du Québec à Trois-Rivières, Trois-Rivières, Canada}
\email{Benjamin.Gauthier@uqtr.ca \qquad Alain.Goupil@uqtr.ca \qquad Fadel.Toure@uqtr.ca}
}
\def\titlerunning{Generation of Maximal Snake Polyominoes  Using a Deep Neural Network}
\def\authorrunning{B. Gauthier, A. Goupil \& F. Toure}

\begin{document}
	\maketitle

\begin{abstract}

    Maximal snake polyominoes  are difficult to study numerically in large rectangles, as computing them requires the complete enumeration of all snakes for a specific rectangle size, which corresponds to a brute force algorithm.
    This hinders the study of maximal snakes in larger rectangles.
    Moreover, most enumerable snakes lie in small rectangles, obscuring large-scale patterns.
    In this paper, we investigate the contribution of a deep neural network to the generation of maximal snake polyominoes  from a data-driven training, where the maximality and adjacency constraints are not encoded explicitly, but learned. To this extent, we experiment with a denoising diffusion model, which we referred as Structured Pixel Space Diffusion (SPS Diffusion). We find that SPS Diffusion generalizes from small rectangles to larger ones, generating valid snakes up to $28\times 28$ squares and producing maximal snake candidates on squares close to the current computational limit. The model is, however, prone to errors such as branching, cycles, or multiple snake components. Overall, the diffusion model is promising and suggests that complex combinatorial objects can be understood by deep neural networks, which is useful in their investigation.
\end{abstract}

\section{Introduction}

    In the past $40$ years polyominoes and polycubes have been the object of intense investigations from different points of view: combinatorial, algorithmic and statistical physics (\cite{guttmann2009history} and ref. therein, \cite{masse2018saturated}). The problem of finding the length of maximal snake polyominoes  in an arbitrary rectangle is still open.
    The most effective algorithmic approach currently available consists in enumerating all snake polyominoes  and then selecting the maximal ones.

    Early foundational work on such brute force algorithms was conducted by Jensen \cite{Jensen2000}, who used transfer-matrix methods and enumeration techniques to count polyominoes and tree-like polyominoes. Later work with Guttmann \cite{Jensen2001} analyzed the asymptotic behavior of polyomino counts using tools from statistical mechanics, showing that the number of polyominoes grows by a factor of about $4.06$ when the number of cells increases by one.
    Since  there exists families of snake polyominoes that are in bijection with Fibonacci numbers (\cite{goupil2018partially}), the number of snakes polyominoes grows exponentially with respect to their size and since there is an injection from the set of snake polyominoes of length $n$ into the set of self-avoiding walks (\cite{goupil2018partially}) the growth constant of snake polyominoes is at most that of self-avoiding walks which is known to be $2.6$ (\cite{jacobsen2016growth}).
    There is thus an inherent difficulty in generating snakes by exhaustive enumeration, especially in large rectangles with size  beyond the range where brute force enumeration is feasible.

    Barequet and Ben-Shachar \cite{GilGill2024} revisited polyomino enumeration by introducing bounding boxes rotated by $45^\circ$ instead of axis-aligned rectangles. This reformulation improved brute force enumeration and allowed polyominoes to be enumerated up to size $70$, compared to Jensen's previous limit of $56$. Nevertheless, the enumeration complexity remains exponential, highlighting the difficulty of counting and generating maximal snakes in large rectangles.

    Another issue is that we observed (without proof) that the general structure of maximal snake polyominoes becomes more apparent in larger rectangles compared to small rectangles. In our brute force computations, we were able to reach sizes of up to $15 \times 15$ for squares, where the general pattern begins to emerge, but significantly larger rectangles remain out of reach of conventional methods.
    It becomes relevant to observe maximal snake polyominoes in large rectangles in the hope of formulating conjectures of exact expressions about their size and their enumeration.

    A theoretical upper bound on the maximal length of snake polyominoes is provided by the 2/3 theorem \cite{blondin2025maximal}. While these upper bounds provide useful theoretical guidance, they do not help overcome the computational difficulties of generating maximal snake polyominoes in large rectangles.

    Given these limitations, deep neural networks (DNNs) offer a possible way to study snakes in larger rectangles.
    Indeed, DNNs have polynomial complexity, which makes them theoretically  faster than the current brute force algorithms used for generating  snake polyominoes  of maximal length. While DNNs have achieved remarkable success in domains such as image generation and natural language processing \cite{vaswaniAttention2023,ho2020denoisingdiffusionprobabilisticmodels}, their application in discrete combinatorial objects and mathematical structure discovery remains relatively unexplored.

    In this work, we generate snake polyominoes using a Structured Pixel Space Diffusion (SPS Diffusion) model. Our goal is not to replace exact methods for enumerating snake polyominoes, but rather to evaluate whether deep neural networks can learn the structural constraints of snake polyominoes  and produce candidates for maximal snakes, thus offering an additional tool to classical enumeration techniques and, as a byproduct, an argument for the existence of combinatorial structure of these objects.
    Section \ref{sec:preliminaries} introduce the necessary concepts. In section  \ref{sec:Methodology for the SPS Diffusion model}
we introduce two neural network and compare them. In section  \ref{sec:Results for the SPS Diffusion model}. we present the results obtained.

    \section{Preliminaries} \label{sec:preliminaries}

    \subsection{Snake polyominoes  definitions}

    We recall the definitions needed below.
    For a rigorous definition of snake polyominoes, see \cite{lheureuxRecherche2022}.

\begin{figure}[ht]
    \centering
    \begin{subfigure}[t]{0.3\textwidth}
        \centering
        \includegraphics[width=0.5\linewidth]{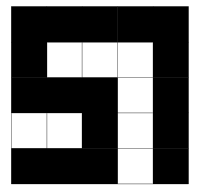}
        \caption{Grid representation of a snake.}
        \label{subfig:multiple_snakes}
    \end{subfigure}
    \begin{subfigure}[t]{0.3\textwidth}
        \centering
        \includegraphics[width=0.5\linewidth]{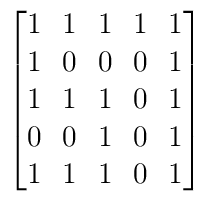}
        \caption{Matrix representation of a snake.}
        \label{subfig:multiple_snakes_tensor}
    \end{subfigure}
    \caption{Different representations of snake polyominoes.}
    \label{fig:rep_polyomino_as_tensor}
    \end{figure}

    A snake polyomino, or more briefly a snake, is a polyomino where all cells are of degree two except for two cells of degree one, called the head and tail. Here, a maximal snake means a snake of maximum area in a given rectangle. In a grid representation, living cells (black) correspond to the snake body and dead cells (white) to empty cells (Fig.~\ref{subfig:multiple_snakes}). We use a matrix representation where living cells are labeled $1$ and dead cells $0$ (Fig.~\ref{subfig:multiple_snakes_tensor}).
    Snakes, and their matrix representations, are contained in discrete rectangles of size $H\times W$.
    From graph theory, the degree of a cell $x$ is the number of edge-adjacent living cells. We assume the standard 4-connexity of square cells, where each cell has at most four neighbors \cite{diestelGraph2025}.

\begin{figure}[ht]
        \centering
        \includegraphics[width=0.6\linewidth]{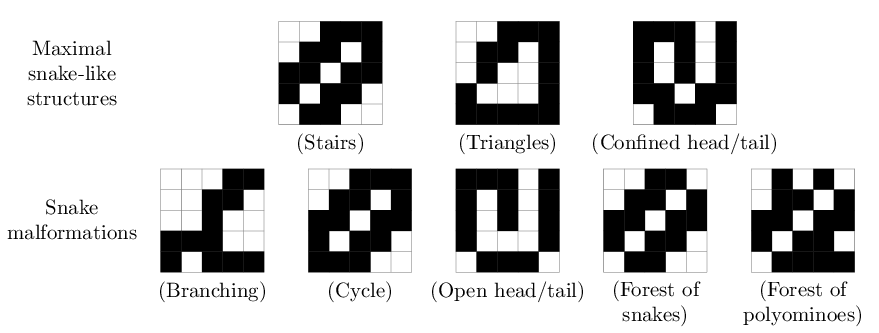}
        \caption{Maximal snake-like structures (above) - stairs, a triangle of dead cells and the confined tail/head. Snake malformations (below) - branching,  cycle,  forest of snakes and an open head/tail. Note that for the open head/tail, the disjointed segments should be considered part of the same snake, as if the rectangle was bigger than shown here.}
        \label{fig:structures}
    \end{figure}

    To qualitatively assess generated snakes, we identify recurring snake-like structures, i.e., patterns commonly observed in maximal snakes (Fig.~\ref{fig:structures}). Examples include stairs, which provide a compact folding of cells of degree at most two \cite{blondin2025maximal}, triangles of six dead cells, and configurations where the head or tail folds compactly. Such structures frequently appear in maximal snakes, although they are not exclusive. For example, head and tail are  frequently observed along the borders and corners of the rectangle.

    Conversely, we identify snake malformations, which are structures incompatible with valid snakes (Fig.~\ref{fig:structures}). Examples include branching, cycles (absence of head and tail), open heads/tails and forests of snakes where several disconnected components appear. The presence of such structures generally indicates that the snake is not maximal. Note that an open head/tail doesn't imply that the snake is not maximal, rather that it may not be.
    It is also useful to know that density is in general greater along the edges of the rectangle than in the center. This fact has been proved for $2D$-words (\cite{blondin2025maximal}).  Border density is thus a good indicator that a snake is maximal.

    \subsection{Deep neural network definitions}

    In the SPS Diffusion model to be  introduced later, data are represented as matrices. The source ($\mathbf{src}$) is the input of the model, the target ($\mathbf{tgt}$) the expected output, and the output ($\mathbf{out}$) the model prediction. The model itself is a parameterized function represented by a deep neural network and learned from data. 

    Training consists in minimizing a loss function that measures the discrepancy between $\mathbf{out}$ and $\mathbf{tgt}$. The parameters of the neural network are updated using gradient descent combined with backpropagation, which computes the gradient of the loss with respect to the model parameters \cite{ruderOverview2017}.

    A fundamental layer of diffusion model is the convolutional layer, which is itself a neural network component that applies small learnable filters (kernels) across an input (such as an image) to extract local patterns. Each filter slides over the input and computes dot products, producing feature maps that detect structures such as edges, textures, or shapes. Convolutional layers are central to diffusion models and enable the efficient capture of spatially localized features, which is particularly important for modeling snakes \cite{ronnebergerU-net2015}.

    Finally, we define self-attention, which is a neural network mechanism that allows a model to contextualize  different parts of an input to one another. Given a set of input representations, the layer computes queries, keys, and values, and determines how strongly each element should attend to the others through similarity scores (attention scores). This produces weighted combinations of the values, enabling the model to capture long-range dependencies and global structure \cite{vaswaniAttention2023}. In spatial settings such as snake polyomino generation, positional information is incorporated using 2D Rotary Positional Embeddings (RoPE-2D) \cite{suRoformer2023}, which encode relative spatial relationships while preserving the ability of attention to model interactions across the rectangle.

    \section{Methodology for the SPS Diffusion model}
    \label{sec:Methodology for the SPS Diffusion model}

    \subsection{Snake generation} \label{sec:diffusion_generation}

    We investigate the use of a denoising diffusion probabilistic model (DDPM) for generating maximal snakes.
    We call our variant \emph{structured pixel space diffusion} (SPS Diffusion), for reasons detailed in Sec.~\ref{sec:diffusion_differences}.

    The process of generating polyominoes with SPS Diffusion is directly analogous to image generation with a standard DDPM and it is described as follows:

    \begin{enumerate}
        \item An image of pure noise, sampled from a normal distribution, is initialized and given to the model.
        \item At each step, the model predicts the amount of noise to be removed from the sample in order to move closer to a structured configuration.
        \item The predicted noise is subtracted from the current sample.
        \item The updated sample is fed back into the model, and the process is repeated for a given number of steps. A larger number of steps typically yields higher-quality outputs.
    \end{enumerate}

    \begin{figure}[b]
        \centering
        \includegraphics[width=.8 \linewidth]{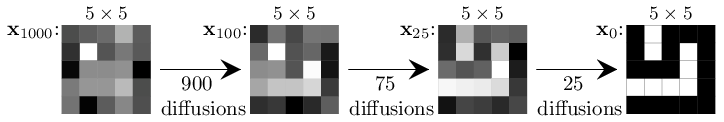}
        \caption{The backward diffusion process to obtain a maximal snake  using an SPS Diffusion model. The arrows indicate the number of diffusion steps used for the next image $(\mathbf{x}_t)$. The initial pure noise image $\mathbf{x}_{1000}$  is diffused into a maximal snake image $\mathbf{x}_0$.}
        \label{fig:gen_steps_diffusion}
    \end{figure}

    Fig.~\ref{fig:gen_steps_diffusion} illustrates the evolution of an image from pure noise to a maximal snake.
    This process is called \textit{backward diffusion}. 
    As in standard DDPMs, the backward diffusion process is described by Eq.~\eqref{eq:backward_diffusion}.

    \begin{equation}
        \mathbf{x}_{t-1} = \frac{1}{\sqrt{\alpha_t}}\left(\mathbf{x}_t - \frac{\beta_t}{\sqrt{1 - \bar \alpha_t}}\epsilon_\theta(\mathbf{x}_t,t)\right) + \sqrt{\frac{1-\bar\alpha_{t-1}}{1 - \bar\alpha_t}\beta_t} \epsilon \label{eq:backward_diffusion}
    \end{equation}

    \noindent where $\mathbf{x}_t$ is the image at step $t$, $\epsilon_\theta(\mathbf{x}_t,t)$ is the prediction of the model and $\epsilon \sim \mathcal{N}(0,1)$ is a statistical variable following a normal distribution of average 0 and variance 1 \cite{ghojoghBenyaminDiffusionSurvey}.
    The image $\mathbf{x}_0$ is the fully denoised image (a maximal snake ideally), and the image $\mathbf{x}_T$ is the pure noise image (noisy polyomino).
    The other symbols are constants from the forward diffusion process; see Sec.~\ref{sec:diffusion_training} and \cite{ghojoghBenyaminDiffusionSurvey}.

    \subsection{Differences between SPS Diffusion and Stable Diffusion} \label{sec:diffusion_differences}

    SPS Diffusion operates similarly to Stable Diffusion: a noisy image is iteratively denoised until the final output emerges. However, maximal snakes are structurally much simpler than natural images. Polyominoes are represented in binary form, with a single feature channel taking values 0 (white) or 1 (black). In practice, this channel initially takes values in $\mathbb{R}$ and is rounded to the nearest integer only at the end, which explains the gray cells visible during the backward diffusion process in Fig.~\ref{fig:gen_steps_diffusion}.

    In contrast, natural images typically use three feature channels (red, green, blue), each ranging from 0 to 255, yielding over 16 million possible colors. Moreover, unlike Stable Diffusion, the model is not conditioned on text or image prompts: the generation always converges to a specific family of images, namely maximal snakes. Finally, snake grids contain far fewer pixels than photographic images; in practice they are smaller than $64\times64$, whereas natural images often exceed $512\times512$ pixels.
    These differences necessitate  modifications to the standard Stable Diffusion architecture.

    \textbf{First}, our model does not include a Contrastive Language–Image Pretraining (CLIP) encoder.
    In Stable Diffusion, CLIP is used to project text and image data into a shared embedding space, enabling training on text–image pairs.
    CLIP can also serve as a loss function to measure how well a generated image matches a textual prompt \cite{radfordLearning2021}.
    Since maximal snake polyomino generation involves no prompts, this module is unnecessary.

    \textbf{Second}, the model is composed solely of a U-Net.
    As in Stable Diffusion, the U-Net predicts the noise to be removed from the current sample.
    To achieve this, the U-Net reduces the spatial resolution of the input as it applies successive convolutional layers, a process known as downsampling \cite{ronnebergerU-net2015}.
    However, because snake polyomino images are already small, excessive downsampling would collapse the structural information in the image entirely.
    While a typical Stable Diffusion U-Net may employ 5 or 6 levels of downsampling, SPS Diffusion is restricted to only 3 levels to preserve structural fidelity.
    For these reasons, this U-Net is better referred to as a \emph{mini U-Net}.
    Each level of the mini U-Net applies several convolutions to the image, followed by a self-attention computation which uses RoPE 2D.

    \textbf{Third}, the model does not require a Variational Autoencoder (VAE).
    In Stable Diffusion (or more generally, Latent Diffusion models), the VAE compresses images from pixel space into a latent space to reduce computational cost, since the U-Net cannot practically process very large images directly \cite{rombachHigh-resolution2022}.
    Snake images, however, are already very small, making such a compression step useless.

    Since there is no VAE, there is no latent space, thus only leaving the pixel space in which the polyomino images live. Also, because this model is specialized in  denoising  2D structures, we call it Structured Pixel Space Diffusion (SPS Diffusion). This model is architecturally simpler than Stable Diffusion, having only a single mini U-Net.

\begin{figure}[ht]
        \centering
        \includegraphics[width=\linewidth]{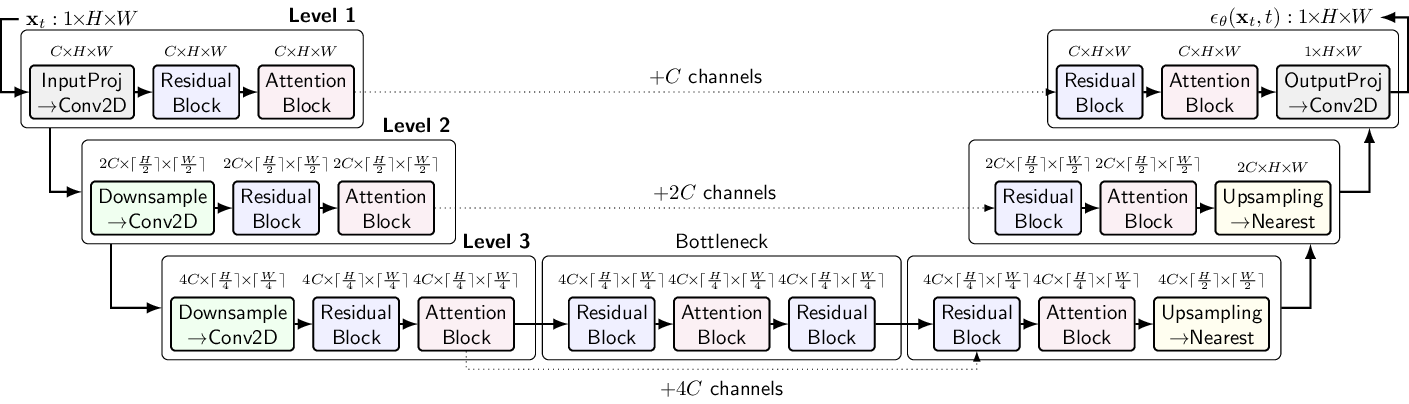}
        \caption{Architecture of the SPS Diffusion model. The model is given $\mathbf{x}_t$, which returns the noise to remove $\epsilon_\theta(\mathbf{x}_t,t)$ (Sec.~\ref{sec:diffusion_training} discusses this further). The labels for the image size indicate the dimensions after the operation of the block below has been applied.}
        \label{fig:sps_diffusion_architecture}
    \end{figure}

    Fig.~\ref{fig:sps_diffusion_architecture} presents the SPS Diffusion model architecture.
    The model is built from residual blocks, which let information pass through while refining features, and attention blocks, which help to focus on the most important parts of the input.
    The network has multiple levels which process the image at different resolutions: lower levels handle fine details, higher levels capture larger patterns.
    The bottleneck is the middle stage where the image is the smallest and most abstract, compressing information before reconstruction.
    Downsampling is done using a specific way of doing convolution, where the kernel is slid across the image 2 pixels at a time, instead of 1 pixel at a time.
    Upsampling is done using nearest-neighbor interpolation, which enlarges the image back toward its original size.
    Skip connections link matching levels in the encoder and decoder, adding the feature channels from the encoder to the decoder to help recover fine details more accurately.

    \subsection{Training} \label{sec:diffusion_training}

    Sec~\ref{sec:diffusion_generation} explains how polyomino images are generated using the SPS Diffusion model. This process is called backward diffusion, because noise is removed from the image instead of being added. When training DDPM models, we use \textit{forward diffusion}, where noise is added to a known solution $\mathbf{x}_0$, a maximal snake, in order to generate the $\mathbf{src}$ data for training. In other words, the noisy polyominoes $\mathbf{x}_t$ ($\mathbf{src}$), are computed in training. Eq.~\eqref{eq:forward_diffusion} shows how $\mathbf{x}_t$ is computed. The constant $\bar \alpha_t$ is related to the mean and variance parameters of the normal distribution governing the forward diffusion process. See \cite{ghojoghBenyaminDiffusionSurvey} for the full definition.

    \begin{equation}
        \mathbf{x}_t = \sqrt{\bar \alpha_t}\mathbf{x}_0 + \sqrt{1 - \bar \alpha_t} \epsilon \label{eq:forward_diffusion}
    \end{equation}

    Since the model estimates the noise term $\epsilon_{\theta}(\mathbf{x}_t, t)$ (see Eq.~\eqref{eq:backward_diffusion}) that should be removed from $\mathbf{x}_t$ to recover a sample closer to $\mathbf{x}_0$,
     the training target $\mathbf{tgt}$ corresponds to the actual noise $\epsilon$ that was added during the forward diffusion step linking $\mathbf{x}_{t-1}$ and $\mathbf{x}_t$ (Eq.~\eqref{eq:forward_diffusion}).
    Thus, the training pair $(\mathbf{src}, \mathbf{tgt})$ consists of $\mathbf{src} = \mathbf{x}_t$, the noisy polyomino at timestep $t$, and $\mathbf{tgt} = \boldsymbol{\epsilon}$, the true noise added at that step.

    The SPS Diffusion model presented here was trained by computing the $\mathbf{src}$ (a noisy polyomino) using known maximal snakes for rectangles of dimension $W \leq 50$ and $H \leq 14$. The noise $\epsilon$ that was added to the maximal snakes during the computation is the $\mathbf{tgt}$, as per Eq.~\eqref{eq:forward_diffusion}.

    For the loss function, we use the Mean Squared Error (MSE), which was the only tested loss to consistently produce organized snake-like structures (Sec.~\ref{sec:preliminaries}) across all rectangle sizes tested. Losses based on structural penalties, such as degree-three cells or disconnected components, were less effective because early generated images are too noisy for these criteria to be reliable. MSE is defined as follows:

    \begin{equation}
        {\rm MSE}({\rm \mathbf{out}},{\rm \mathbf{tgt}}) = \frac{1}{N}\sum^N_{ij} ({\rm \mathbf{tgt}}_{ij}-{\rm \mathbf{out}}_{ij})^2
    \end{equation}

\section{Results for the SPS Diffusion model}
 \label{sec:Results for the SPS Diffusion model}

\begin{figure}[t]
        \centering
        \includegraphics[width=0.8\linewidth]{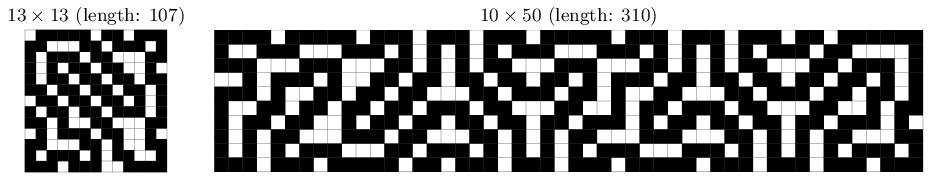}
        \caption{Two maximal snakes generated by the SPS Diffusion model in rectangles seen during training.}
        \label{fig:diffusion_seen}
    \end{figure}

     We evaluate SPS Diffusion on different rectangle sizes.
    We consider a SPS Diffusion model trained on over 200 different rectangle sizes, ranging from small $4\times4$ squares to $50\times10$ rectangles.
    We start with rectangles for which the sizes have been seen by the model during training, as shown in Fig.~\ref{fig:diffusion_seen}.
    SPS Diffusion has little difficulty generating maximal snake polyominoes  in all rectangle sizes seen during training.

\begin{figure}[t]
        \centering
        \includegraphics{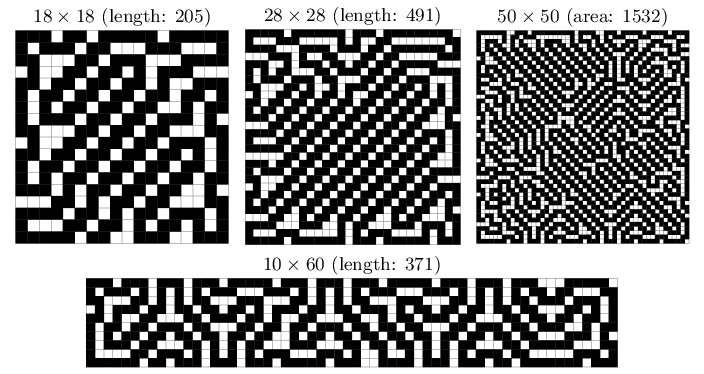}
        \caption{Results generated by the SPS Diffusion model in squares not seen during training.}
        \label{fig:diffusion_not_seen}
    \end{figure}

    Looking at rectangles not seen by the model during training in Fig.~\ref{fig:diffusion_not_seen}, SPS Diffusion is able to generalize its understanding of snake polyominoes  in small rectangles to much larger sizes at inference time, that is, generate snake-like structures even in large rectangles where the positional encoding is technically different.
    Furthermore, we were able to generate snake polyominoes in squares for sizes up to $28\times28$, which is bigger than the maximum size used in training for squares ($14\times14$). 
    Similarly, for many rectangles, ranging from near squares to very elongated rectangles of size $10 \times 60$, we were also able to generate several snakes.
    When trying to generate snakes in very large rectangle sizes, we hit the limitations of the model. Fig.~\ref{fig:diffusion_not_seen} shows a forest of polyominoes generated in a large square of size $50\times50$.
    Although we observe snake-like structures in the square, the structure itself is not a  snake. In some cases, we obtain a forest of  snakes, in others, we obtain forests of polyominoes, which are not snakes.
    Considering the complexity of maximal snakes, the fact that SPS Diffusion is able to get close to a valid solution is encouraging.

    Moreover, we were able to generate maximal snakes for rectangle sizes up to $17\times 17$, which is the current computational limit for squares (\cite{oeis-a331968}).
    After this bound, we cannot know for certain if the generated snakes are maximal.
    However, we can estimate a lower bound on the minimal length of maximal snakes and, for several rectangle sizes, we can sometimes exceed the lower bound by up to $3$ cells.
    We observe that as the rectangle size increases, it becomes harder for the model to generate a snake, and even harder for a maximal snake. In other words, the error rate increases with the rectangle size. We also observe that the SPS Diffusion model makes fewer errors on elongated rectangles than on squares, since the corresponding snakes are more predictable and therefore easier to generate.

    The aformentioned limit of $28 \times 28$ for squares containing a snake is a consequence of this fact.
    For example, for more than  $25\; 000$ images generated, no snakes were produced for a $29\times29$ square, only forests of snakes and forests of polyominoes.
    We expect that with  more calculation time, snakes could be generated for larger square sizes, but then we would encounter the initial limitations of the brute force algorithm, that is, calculation time.

    While the SPS Diffusion model generates valid snakes in large rectangles, it has also produced snakes whose lengths exceed the best previously reported lower bounds. These bounds are rigorous rather than conjectural, as snakes of these lengths have been explicitly constructed, often by hand.

    Before these results were obtained, the SPS Diffusion model was trained on a smaller dataset of rectangles smaller than $20\times10$. In this setting, the model generated snakes in squares up to $23\times23$. After adding squares up to $14\times14$ (and other similarly sized rectangles) to the training data, the model generated snakes up to $28\times28$. This improvement suggests that the model benefits substantially from additional data. Consequently, one possible approach is to generate new snakes with the model itself and include them in the training set for retraining. Such an iterative procedure could improve predictions in larger rectangles, although it cannot guarantee that the produced snakes are maximal and mainly serves to guide theoretical investigation of maximal snake lengths.

    This procedure can in principle be repeated to generate increasingly large snakes, although it currently requires significant computational time. The SPS Diffusion model also makes several errors and produces forbidden patterns such as cells of degree $3$ (branching), cycles, forests of snakes and polyominoes, requiring more computation time in order to generate valid snakes, which, even then, yields only a few snakes, and even less maximal snakes.

    Moreover, in larger rectangles there is no guarantee that the positional encoding used for smaller rectangles (RoPE 2D in the attention layers) will remain effective. The current mini-U-Net architecture is also fixed, with only three downsampling layers, which may be insufficient to capture global positional context in larger rectangles. This limitation could be mitigated by training a larger U-Net while excluding smaller rectangles from the training data. Overall, further optimization of the model implementation and evaluation on larger rectangles will be necessary for long-term applicability.

    \section{Conclusion}

    In this work, we explored whether a deep neural network can generate maximal snake polyominoes beyond the range of brute force enumeration.
    We used a diffusion model called Structured Pixel Space Diffusion (SPS Diffusion).

    Our experiments showed that SPS Diffusion generates good results,
    as it can easily reproduce the data it was trained with, as well as generalizing to larger rectangles.
    More specifically, for squares, snakes were generated up to $28\times28$, and their length was close to the previously known lower bound. For some configurations, the previously known bound was exceeded by one cell, which demonstrates the potential of the approach when given more training and refinements.
    Its main weakness is the generation of malformations, including branching, cycles, and forests.
    Also the likelihood of generating a valid snake decreases as the rectangle size increases. 
    Overall, the results suggest that diffusion models can learn nontrivial combinatorial structures from examples, even when constraints are not explicitly encoded.

    Several directions remain open. On the modeling side, reducing invalid generations and improving scalability to larger rectangles are key objectives, for instance by modifying the architecture and the positional encoding of the cells. Alternatively, the loss function could also be tweaked to speed up and improve the training.

    On the combinatorial side, a more systematic analysis of recurring substructures in the generated candidates could help guide future conjectures about maximal snake size in large rectangles. Overall, deep neural network models, and SPS Diffusion in particular, provide a new numerical tool to explore maximal  snakes beyond the range of exhaustive enumeration.

\bibliographystyle{eptcs}
\bibliography{references}

\end{document}